\documentclass[reqno]{article}
\usepackage{tikz}
\usepackage{pgf}
\tikzset{%
    vtx/.style={draw,circle,ultra thin,white,fill=black,inner sep=1.5pt,text width=2mm},
    bigvtx/.style={draw,circle,ultra thin,white,fill=black,inner sep=2pt,text width=2.5mm},
    edg/.style={ultra thick},
    spedge/.style={rounded corners,line width=5pt,red!70}, 
    thedge/.style={} 
}
\tikzset{/pgf/foreach/parse=true}
\pgfdeclarelayer{back}
\pgfdeclarelayer{fore}
\pgfsetlayers{back,main,fore}
\usepackage{amsmath,amsthm,amssymb}
\usepackage{graphicx}
\graphicspath{ {./Steiner/} }
\newtheorem{theorem}{Theorem}[section]
\newtheorem{lemma}[theorem]{Lemma}
\newtheorem{proposition}[theorem]{Proposition}

\newtheorem{conjecture}[theorem]{Conjecture}

\newtheorem{claim}[theorem]{Claim}

\def\beq{\begin{equation}}\def\eeq{\end{equation}}
\def\beqn{\begin{eqnarray}}\def\eeqn{\end{eqnarray}}

\def\qed{\ifhmode\unskip\nobreak\fi\quad\ifmmode\Box\else$\Box$\fi}

\begin{document}
\title{Proper edge colorings of planar graphs with rainbow $C_4$-s}
\author{Andr\'as Gy\'arf\'as\thanks{Alfr\'ed R\'enyi Institute of Mathematics, Budapest, P.O. Box 127, Budapest, Hungary, H-1364.
Corresponding author: \texttt{gsarkozy@cs.wpi.edu}} \thanks{Research supported in part by
NKFIH Grant No. K132696.} \and Ryan R. Martin\thanks{Iowa State University, Ames, Iowa, USA. Research partially supported by Simons Foundation Collaboration Grant for Mathematicians \#709641 and this work was done in part while this author was a Distinguished Guest Fellow of the Hungarian Academy of Sciences at the R\'enyi Institute for Mathematics.} \and Mikl\'os Ruszink\'o\footnotemark[1] \footnotemark[2]
\thanks{Department of Sciences and Engineering, Sorbonne University Abu Dhabi; Faculty of Information Technology and Bionics, P\'azm\'any P\'eter Catholic University.} \and G\'{a}bor N. S\'ark\"ozy\footnotemark[1] \footnotemark[2]
\thanks{Computer Science Department, Worcester Polytechnic Institute, Worcester, MA, USA, Research supported in part by
NKFIH Grant No. K117879.}}

\maketitle
\date{}
\begin{abstract}

We call a proper edge coloring of a graph $G$ a B-coloring if every 4-cycle of $G$ is colored with four different colors. Let $q_B(G)$ denote the smallest number of colors needed for a B-coloring of $G$. Motivated by  earlier papers on B-colorings, here we consider $q_B(G)$ for planar and outerplanar graphs in terms of the maximum degree $\Delta = \Delta(G)$. We prove that
$q_B(G)\le 2\Delta+8$ for planar graphs, $q_B(G)\le 2\Delta$ for bipartite planar graphs and $q_B(G)\le \Delta+1$ for outerplanar graphs with $\Delta \ge 4$.
We conjecture that, for $\Delta$ sufficiently large, $q_B(G)\le 2\Delta(G)$ for planar $G$ and $q_B(G)\le \Delta(G)$ for outerplanar $G$.\\
{\bf Keywords:} planar graphs, graph colorings.

\end{abstract}

\section{Introduction}

A graph $G$ is called {\em planar} if it can be drawn on the Euclidean plane so that any two edges intersect only at their ends. A special class of planar graphs are {\em outerplanar} graphs which admit a plane embedding such that all of its vertices lie on the boundary of the same face. A $2$-connected outerplanar graph $G$ is a Hamiltonian cycle with some chords. $K_n$ denotes the complete graph on $n$ vertices and $K_{m,n}$ denotes the complete bipartite graph between two vertex sets of sizes $m$ and $n$. Denote by $K_n - e$ the graph we get if we remove the edge $e$ from $K_n$.
The {\em Cartesian product} of two graphs $G$ and $H$, denoted by $G\square H$, is a graph with vertex
set $V(G)\times V (H)$, and $((a,x),(b,y))\in E(G\square H)$ if either $(a,b)\in E(G)$ and $x = y$, or $(x,y)\in E(H)$
and $a = b$. For a graph $G$ and a subset
$U$ of its vertices, $G|_U$ is the {\em restriction} of $G$ to $U$. $N(v)$
is the set of neighbors of $v\in V$. Hence $|N(v)|=deg(v)=deg_G(v)$,
the {\em degree} of $v$.

A graph $G$ has a graph $H$ as a {\em minor} if $H$ can be obtained from a subgraph of $G$ by contracting edges. Furthermore, $G$ is called {\em $H$-minor free} if $G$ does not have $H$ as a minor. A graph $G$ has a graph $H$ as a {\em subdivision} if $H$ can be obtained from a subgraph of $G$ by replacing internally-disjoint induced paths with edges. It is well-known that a graph $G$ is planar if and only if $G$ is $K_5$-minor free and is $K_{3,3}$-minor free (Wagner's theorem \cite{WAG}) and that $G$ is planar if and only if $G$ has no subdivision of $K_5$ or of $K_{3,3}$ (Kuratowski's theorem \cite{KUR}).
It was shown in \cite{CH} that a graph $G$ is outerplanar if and only if $G$ is $K_4$-minor free and $K_{2,3}$-minor free. It is an easy exercise to use Kuratowski's theorem to prove that $G$ is outerplanar if and only if $K$ has no subdivision of $K_4$ or of $K_{2,3}$. A regular graph of degree three is called {\em cubic}; a graph with maximum degree at most three is called {\em subcubic}.

An edge coloring of a graph $G$ is {\em proper} if incident edges of $G$ must receive different colors. Recently two of the authors studied proper edge colorings with the additional requirement that  every 4-cycle of $G$ must be {\em rainbow}, i.e., colored with four different colors. We called these colorings {\em B-colorings} and defined $q_B(G)$ as the smallest number of colors needed for a B-coloring of $G$ \cite{GYS1}. The motivation to study $q_B(G)$ comes partly from earlier works \cite{FGLS,FGS} and partly from its relation to the famous $(7,4)$-problem of Brown, Erd\H{o}s and S\'os \cite{BETS}. We asked in \cite{GYS1} whether any graph $G$ with $n$ vertices and with $q_B(G)=cn$ has $o(n^2)$ edges. We showed that a positive answer to this question would imply a positive answer to the $(7,4)$-problem as well: any triple system on $n$ points with no 4 triples on 7 vertices has $o(n^2)$ triples \cite{BETS}.
It is worth noting that B-colorings are similar to the well studied concept of {\em star-edge colorings}: proper edge colorings where the union of any two color classes does not contain paths or cycles with four edges. These colorings were defined in \cite{DMS} (though having appeared already in \cite{EGY}). It was shown in \cite{GYS1} that (similarly to B-colorings) the following statement would give a positive answer to the $(7,4)$-problem. Any graph $G$ with $n$ vertices and with star-edge colorings with $cn$ colors has $o(n^2)$ edges.

Here we study $q_B(G)$ for planar graphs. Note that star-edge colorings of planar graphs are  well-studied and provided many open problems (see e.g. \cite{BLMSS,DYZC,WWW}, or more results in the survey  \cite{LS}).

\subsection{The extended line-graph $L^+(G)$}

The {\em chromatic number} of a graph $G$ is denoted by $\chi(G)$ (the minimum number of colors needed in a proper vertex-coloring) and the {\em clique number}
is denoted by $\omega(G)$ (the size of the largest clique).
The parameter $q_B(G)$ can be naturally translated to finding the chromatic number of the graph that we call $L^+(G)$, the ``extended'' line-graph of $G$. Traditionally, $L(G)$, the {\em line-graph} of $G$ is defined by representing the  edges of $G$ by vertices and defining two vertices adjacent if and only if the corresponding edges are incident in $G$. In the {\em extended line graph}, $L^+(G)$ we add further edges to $L(G)$ as follows: we add an edge between two non-adjacent vertices of $L(G)$ if the corresponding (independent) edges in $G$ are in a 4-cycle of $G$. We denote by $F(G)$ the graph induced by this set of new edges. As a result, $L^+(G)=L(G)\cup F(G)$.
\begin{proposition}\label{trans} For any graph $G$, $q_B(G)=\chi(L^+(G))$.
\end{proposition}

Proposition \ref{trans} leads to  \beq\label{inequal} \omega(L^+(G))\le \chi(L^+(G))=q_B(G),\eeq
providing a natural lower bound of $q_B(G)$.
In the cases of planar and outerplanar graphs $G$ we have some restrictions on $F(G)$. These are stated as two lemmas needed for the proofs of our main results.

\begin{lemma}\label{planarF}Let $G$ be a planar graph.
\begin{itemize}
\item [(i)] If $F(G)$ has a subgraph $H$ with $k\geq 2$ vertices $x_1,\dots,x_k$ such that the corresponding $k$ edges in $G$ are vertex disjoint (i.e. they form a matching),
then $G$ contains the subgraph $H\square K_2$ in which each $x_i$ corresponds to an edge $e_i$ between the two copies of $H$, called connector edges.
\item [(ii)] $\omega(F(G))\le 3$.
\end{itemize}
\end{lemma}

\begin{lemma}\label{outerplanarF} If $G$ is outerplanar then  $\omega(F(G))\le 2$.
\end{lemma}

\subsection{B-coloring of planar and outerplanar graphs}

Our main results are the estimates of $q_B(G)$ for planar and outerplanar graphs.

\begin{theorem}\label{planar}
Let $G$ be a planar graph with maximum degree at most $\Delta$. Then we have
\begin{itemize}
\item[(i)] $q_B(G) \leq 2\Delta + 8$,
\item[(ii)] $q_B(G) \leq 2 \Delta$ if $G$ is bipartite,
\item[(iii)] $\omega(L^+(G)) \leq 2 \Delta$ if $\Delta \ge 6$ and this is sharp.

\end{itemize}
\end{theorem}

Note that (i) and (ii) are also sharp apart from the constant term.
The sharpness of these results is shown by $K_{2,\Delta}$ since all edges must get a different color in any B-coloring.
It is an interesting problem to determine here the true value of $q_B(G)$ in general.
It is tempting to conjecture that $q_B(G) \leq 2 \Delta$ for all planar graphs. However, this is not the case, as evidenced
by the graph $G$ that we get if we remove a perfect matching from $K_6$. Here $\Delta(G)=4$ and in any B-coloring of $G$, all edges must get a different color, thus $q_B(G) = 12 = 2 \Delta + 4$.

For planar graphs with maximum degree at most 3 (subcubic graphs) we can prove the $2\Delta$ bound.

\begin{theorem}\label{subcubic}
If $G$ is a subcubic planar graph, then we have $q_B(G) \leq 6$ and this is sharp.
\end{theorem}

As it is usually the case, for outerplanar graphs we can show stronger results. We state our results for $2$-connected outerplanar graphs since, at cut-vertices, the colorings can be adjusted.

\begin{theorem}\label{outerplanar}
Let $G$ be a $2$-connected outerplanar graph with maximum degree at most $\Delta$. Then we have
\begin{itemize}
\item[(i)] $q_B(C_4)=4, q_B(K_4-e)=5$. If $G$ is distinct from $C_4$ and $K_4-e$, then $q_B(G) \leq \Delta + 1$ and this is sharp for $\Delta=2,4$.
\item[(ii)] $\omega(L^+(G))\le \Delta(G)$ if $\Delta \geq 5$ and this is sharp.
\end{itemize}
\end{theorem}

The sharpness in (i) for $\Delta(G)=2$ is shown by odd cycles; for $\Delta(G)=4$  by a 4-vertex path and a vertex adjacent to all vertices of this path.

\begin{conjecture}\label{conj} For large enough $\Delta= \Delta(G)$ we have
$q_B(G) \leq 2 \Delta$ for any planar $G$ and $q_B(G) \leq  \Delta$ for any outerplanar $G$.

\end{conjecture}

Theorem \ref{planar} part (iii) and Theorem \ref{outerplanar} part (ii) show that Conjecture \ref{conj} is true if $q_B(G)$ is replaced by its lower bound $\omega(L^+(G))$ (see (\ref{inequal})).

\section{Proofs}

$V(G)$ and $E(G)$ denote the vertex-set and the edge-set of the
graph $G$. For a graph $G$ and a subset
$U$ of its vertices, $G|_U$ is the restriction of $G$ to $U$. $N(v)$
is the set of neighbors of $v\in V$. Hence $|N(v)|=\deg(v)=\deg_G(v)$,
the degree of $v$.

\subsection{ Proof of Lemmas \ref{planarF}, \ref{outerplanarF}}

To prove part (i) of Lemma \ref{planarF}, we will use induction on $|V(H)|$.
To launch the induction, note that if $H=K_2$, then by definition the two (vertex disjoint) edges in $G$ corresponding
to the two vertices of $H$ are on a 4-cycle, thus indeed there is a $C_4=K_2\square K_2$ in $G$.
In general, consider a subgraph  $H$ of $F(G)$ with $k$ vertices $x_1,\dots,x_k$ such that the corresponding
$k$ edges in $G$ are vertex disjoint. We may assume that $H$ is connected, otherwise, we can consider its connected components separately.
Select a vertex (say $x_k$) in $H$ that is {\em not} a cut-vertex (it is easy to see that this
must exist).
Let $H_1$ denote the graph we get if we delete $x_k$ from $H$. Note that $H_1$ is still connected, since $x_k$ is not a cut-vertex.

Applying induction, we get a subgraph $H_1\square K_2$ in $G$, where the two $H_1$ copies are joined by the edges $(u_i,v_i), 1\leq i \leq k-1$ in $G$.
Wlog assume that the neighbors of $x_k$ in $H$ are $x_1, \ldots , x_l$ for some $l\geq 1$.
There is a 4-cycle in $G$ containing $x_1=(u_1,v_1)$ and $x_k$. Let $u_k$ be the endpoint of $x_k$ that is connected to $u_1$
on this 4-cycle and let $v_k$ be the other endpoint (it is connected to $v_1$ on the 4-cycle).
Note that by assumption the edge $(u_k,v_k)$ is vertex disjoint from the edges $(u_i,v_i), 1\leq i \leq k-1$.
If $l=1$,
then we are immediately done; we can add $u_k$ to the first copy of $H_1$ (on the vertices $u_i$) and $v_k$ to the second copy.
Assume $l\geq 2$. We claim that all the 4-cycles containing $x_k$ and $x_i, 1\leq i \leq l$ (these must exist)
visit $(u_k,v_k)$ in the same order, i.e. they are $(u_i,v_i,v_k,u_k)$. Then indeed we could add $u_k$ to the first copy of $H_1$ and $v_k$ to the second one.

Assume indirectly, that there is a ``crossing", i.e. say we have the 4-cycles $(u_1,v_1,v_k,u_k)$ but we have $(u_l,v_l,u_k,v_k)$.
Since $H_1$ is connected, there is a path in $H_1$ connecting $x_1$ and $x_l$, wlog assume that this is the path $(x_1, \ldots , x_l)$.
Then $(u_1, \ldots , u_l)$ is a path in the first copy of $H_1$ and $(v_1, \ldots , v_l)$ is a path in the second copy.
In this case the graph $G$ cannot be planar, a contradiction. We can convince ourselves of this fact by finding a subdivision of a $K_{3,3}$ with base points
$\{u_l, v_1, u_k\}$ and $\{v_l,u_1,v_k\}$ (see Figure 1 in the special case $l=3$ and $k=4$).


\begin{figure}[ht]
\centering
\begin{tikzpicture}[scale=1]
	\useasboundingbox (-0.5,-0.75) rectangle (6.5,2.75);
	\coordinate (v3) at (0,0);
	\coordinate (v2) at (2,0);
	\coordinate (v1) at (4,0);
	\coordinate (v4) at (6,0);
	\coordinate (u3) at (0,2);
	\coordinate (u2) at (2,2);
	\coordinate (u1) at (4,2);
	\coordinate (u4) at (6,2);
	\begin{pgfonlayer}{main}
		\draw (u3) node[bigvtx,fill=red,label={[label distance=0pt]above:{\large $u_3$}}]{};
		\draw (u2) node[vtx,fill=black,label={[label distance=2pt]above:{\large $u_2$}}]{};
		\draw (u1) node[bigvtx,fill=blue,label={[label distance=0pt]above:{\large $u_1$}}]{};
		\draw (u4) node[bigvtx,fill=red,label={[label distance=0pt]above:{\large $u_4$}}]{};
		\draw (v3) node[bigvtx,fill=blue,label={[label distance=0pt]below:{\large $v_3$}}]{};
		\draw (v2) node[vtx,fill=black,label={[label distance=2pt]below:{\large $v_2$}}]{};
		\draw (v1) node[bigvtx,fill=red,label={[label distance=0pt]below:{\large $v_1$}}]{};
		\draw (v4) node[bigvtx,fill=blue,label={[label distance=0pt]below:{\large $v_4$}}]{};
	\end{pgfonlayer}
	\begin{pgfonlayer}{back}
		\draw[edg] (u3) -- (u2) -- (u1) -- (u4);
		\draw[edg] (v3) -- (v2) -- (v1) -- (v4);
		\draw[edg] (u3) -- (v4);
		\draw[edg] (v3) -- (u4);
		\foreach\i in {3,2,1,4}{
			\draw[edg] (u\i) -- (v\i);}
	\end{pgfonlayer}	
\end{tikzpicture}
\caption{Finding a subdivision of a $K_{3,3}$ if $u_4$ is connected to $v_3$ and not $u_3$.}
\end{figure}
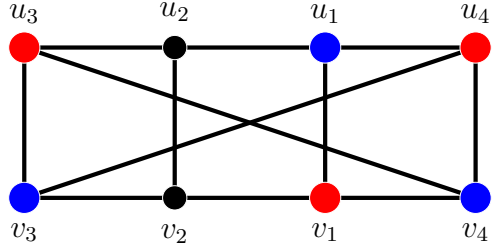

Therefore we have the two copies of $H$, on vertex sets $\{u_i \; | \; 1\leq i \leq k\}$ and $\{v_i \; | \; 1\leq i \leq k\}$ joined by the edges $(u_i,v_i)$ in $G$, forming a subgraph $H\square K_2$ in $G$ as desired.

To prove part (ii), assume indirectly that there is a $K_4$ in $F(G)$ with vertices $w_1,w_2,w_3,w_4$.
Note that in this case the corresponding four edges in $G$ are indeed vertex disjoint, i.e. they form a matching in $G$, so we can apply part (i).
Applying (i) for $H=K_4$, we get a subgraph $H'=K_4\square K_2$ in $G$. However, $H'$ contains a subdivision of $K_5$, contradicting the planarity of $G$.  \qed

To prove Lemma \ref{outerplanarF}, assume indirectly that there is a $K_3$ in $F(G)$ with vertices $w_1,w_2,w_3$. Applying (i) from Lemma \ref{planarF} for $H=K_3$, we get a subgraph $H'=C_3\square K_2$ in $G$, containing a subdivision of $K_4$, contradicting the fact that $G$ is outerplanar. \qed

\subsection{ Proof of Theorem \ref{planar}}

\textbf{Proof of (i):} %
We use induction on the number of vertices and prove that $2\Delta+8$ colors suffice for a proper coloring of a planar graph with maximum degree at most $\Delta$ such that the edges of every $C_4$ are colored with distinct colors. Launching the induction is trivial, we can start from $n=2$. Assume we know the statement for planar graphs of maximum degree at most $\Delta$ and with less than $n$ vertices for some $n\ge 3$.

Let $G$ be a planar graph with $n\ge 3$ vertices and with maximum degree at most $\Delta$.
Since $e(G)\leq 3n-6$, there is a vertex $v$ of degree at most $5$ in $G$. Set $W=N(v)$ and let $H$ be the graph obtained from $G$ by removing $v$ and the edges incident to it. The  graph $H[W]$ has at most seven edges. Indeed, this is trivially true for $|W|\le 4$ and for $|W|=5$ it follows from the planarity of $G$, since $H[W]$ has neither $K_4$ nor $K_{2,3}$ subgraphs.
By the inductive hypothesis, the graph $H$ has the required coloring with at most 2$\Delta+8$ colors. We have to extend this coloring to the (at most 5) edges $\{v,w\}, w\in W$
in such a way that we still have the required coloring.
Consider an order of the vertices $W=\{w_1, w_2, \ldots, w_{|W|}\}$, where we list the degree-two (in $H[W]$) vertices first (if there are any). We will extend the coloring to the edges $\{v,w_1\}$, $\{v,w_2\}$, etc. one-by-one.
Consider a $w_i\in W, 1\leq i\leq |W|$ and set $A_{w_i}=W\setminus \{w_i\}$ and $X_{w_i}=N(w_i)\setminus (A_{w_i}\cup \{v\})$.
Assume that the edges $\{v,w_j\}, 1\leq j <i$ are already colored.
The colors on the following edges will be forbidden to be used for $\{v,w_i\}$ to guarantee the required coloring:

\begin{itemize}
	\item[(a)] 	the edges incident to $w_i$ (at most $\Delta -1$ edges),
	\item[(b)] 	edges within $A_{w_i}$ that are in a four-cycle with the edge $\{v,w_i\}$,
    	\item[(c)]  edges $\{v,w_j\}, 1\leq j <i$ (the edges already colored in the extension, $i-1$ edges).
	\item[(d)] 	edges in the bipartite graph $[A_{w_i},X_{w_i}]$, there are at most $|A_{w_i}|+|X_{w_i}|-1\le \Delta+2$  of these because the graph has no cycle. Indeed, otherwise consider a cycle
$$C=(x_1,w',x_2,\dots,w''),$$ where $x_1,x_2\in X_{w_i}, w',w''\in A_{w_i}$. Observe that the bipartite subgraph of $G$  between $\{v,x_1,x_2\}$ and $\{w_i,w',w''\}$ has eight edges and the (possibly) non-adjacent vertex pair $\{x_2,w''\}$ is connected by a path along $C$, defining a topological $K_{3,3}$, a contradiction.
\end{itemize}

If $w_i$ has degree zero in $H[W]$ then the colors on edges within $A_{w_i}$ are not in a four-cycle containing $\{v,w_i\}$, thus in (b) no edges need to be forbidden, (d) contributes $\Delta+2$, altogether at most $\Delta-1+(i-1)+\Delta+2=2\Delta+i$ colors are forbidden on $\{v,w_i\}$.

If $w_i$ has degree one in $H[W]$ then at most three edges within $A_{w_i}$ are in a four-cycle containing $\{v,w_i\}$, thus (b) contributes at most $3$ and (d) contributes at most $\Delta+1$ (because $|X_{w_i}|\le \Delta-2$). Altogether at most $\Delta-1+3+(i-1)+\Delta+1=2\Delta+i+2$ colors are forbidden on $\{v,w_i\}$.

If $w_i$ has degree two in $H[W]$ then (b) contributes at most $5$ (since $H[W]$ has at most seven edges) and (d) contributes at most $\Delta$ (because $|X_{w_i}|\le \Delta-3$). Altogether at most $\Delta-1+5+(i-1)+\Delta=2\Delta+i+3$ colors are forbidden on $\{v,w_i\}$. This is at most $2\Delta+7$ except if $i=5$. However, in this case all vertices are of degree two in $H[W]$ (since we list the degree-two 
vertices first), and then (b) contributes only 3 (since $H[W]$ has five edges in this case) so similarly as in the previous argument we get at most $\Delta-1+3+(i-1)+\Delta=2\Delta+i+1$ forbidden colors.

If $w_i$ has degree at least three in $H[W]$ then (b) contributes at most $4$ (since $H[W]$ has at most seven edges) and (d) contributes at most $\Delta-1$ (because $|X_{w_i}|\le \Delta-4$). Altogether at most $\Delta-1+4+(i-1)+\Delta-1=2\Delta+i+1$ colors are forbidden on $\{v,w\}$.

Since $i\leq |W|\leq 5$, in all cases, we have at most $2\Delta+7$ forbidden colors for $\{v,w_i\}$. This allows an extension to the required coloring of $G$ with at most $2\Delta+8$ colors 
and concludes the proof.  \qed \\


\textbf{Proof of (ii):} This statement is proved with the technique applied in \cite{FGST} (Theorem 10 in \cite{FGST}).
Let $G$ be a bipartite planar graph with maximum degree at most $\Delta$.
First, by K\"onig's classical theorem we decompose $E(G)$ into $\Delta$  matchings. Then we will color each matching with two colors so that different matchings are colored with distinct colors and there is no four-cycle with two edges of the same color in these matchings. This ensures that the coloring is a B-coloring.

Consider a matching $M$. Denote by $F(M)$ the subgraph of $F(G)$ induced by the edges in $M$.
To get (ii), we need to show that $F(M)$ is bipartite for a bipartite $G$. This follows from Lemma \ref{planarF} (i). Indeed, an odd cycle in $F(M)$
would result in an odd cycle in $G$, a contradiction.  Thus the vertices of $F(M)$ has a proper $2$-coloring, providing a $2$-coloring on $M$. Repeating this for the $\Delta$ matchings (using different color pairs) we get a B-coloring with $2\Delta$ colors.  \qed \\

\textbf{Proof of (iii):} Here we want to prove $\omega(L^+(G)) \leq 2 \Delta$. Let $K$ be a complete subgraph in $L^+(G)$ and let $H$ be the subgraph of $G$ determined by the edge set corresponding to $K$. Let $S$ be the largest complete subgraph of $K$ formed only by edges from $F(G)$. By Lemma \ref{planarF} (ii), $|S|\le 3$.

{\bf Case 1:} $|S|=1$. Here a single edge $e=(u_1,v_1)$ corresponds to $S$ and all further edges of $H$ must meet $e$ thus $\omega(L^+(G))=|E(H)|\le 2\Delta(H)\le 2\Delta$.

{\bf Case 2:} $|S|=2$. Here edges $e_1=(u_1,v_1),e_2=(u_2,v_2)$ correspond to $S$ and $G$  has a 4-cycle $C=(u_1,u_2,v_2,v_1)$ with edges $e_1, e_2$ on it. Clearly no edge $e \in E(H)$  exists with $|e\cap C|=0$ otherwise $|S|>2$. Let $t$ denote the number of edges in $H$ with exactly one vertex on $C$. Since at most 6 edges of $H$ are inside $C$, $|E(H)|\le 6+t$. Suppose indirectly that $|E(H)|>2\Delta$, then we get
$t>2\Delta-6$.

If $e_1,e_2$ both have a vertex, say $x,y$ sending no edges of $H$ to $V(G)\setminus V(C)$, then apart possibly from the edge $e=(x,y)$ all edges of $H$ intersect $T=V(C)\setminus \{x,y\}$. We may assume that $e$ exists, otherwise we have
$|E(H)|\leq 2\Delta$, a contradiction. From the above, we know that more than $2\Delta-6\geq \Delta$ (using $\Delta \geq 6$)  edges of $H$ have exactly one vertex in $T$. Both vertices of $T$ are incident to at least one of these edges, otherwise one of them would be incident to more than $\Delta$ edges of $H$, a contradiction. Moreover, at least one vertex of $T$ is incident to at least two of these edges, thus we can get two independent edges $f,g\in E(H)$ incident to $T$ and meeting $C$ in one vertex. However, now $e,f,g$ are three independent edges of $H$, contradicting the assumption $|S|=2$.

Now assume that one of $e_1,e_2$, say $e_1$ sends from both endpoints at least one edge of $H$  to $V(G)\setminus V(C)$. We get a contradiction with the assumption $|S|=2$, unless $u_1,v_1$ both send just one edge of $H$ to the same vertex $z\in V(G)\setminus V(C)$. If $e_2$ also sends from both endpoints at least one edge of $H$  to $V(G)\setminus V(C)$ then, as before, we have just one further edge from both endpoints of $e_2$, thus $H$ has at most $10\leq 2\Delta$ edges, a contradiction. Otherwise one endpoint of $e_2$, say $u_2$ sends no edges of $H$ to $V(G)\setminus V(C)$. Now we have in $H$ at most $\Delta+5\leq 2\Delta$ edges (using $\Delta\ge 6$), finishing Case 2.

{\bf Case 3:} $|S|=3$. We will show that in this case $|E(H)|\leq 12 \leq 2\Delta$.
By Lemma \ref{planarF} (i) $G$ has the subgraph $G'=C_3\square K_2$. See Figure 2 for a drawing; note that
by Whitney's theorem \cite{WH} since this graph is 3-connected, this drawing is unique. Denote the vertices of $G'$ by $u_i,v_i$ for $i\in [3]$ where the edges $e_i=(u_i,v_i)$ correspond to the vertices of $S$.


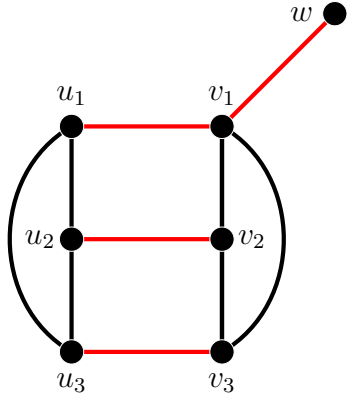
\begin{figure}[ht]
\centering
\begin{tikzpicture}[scale=1]
	\useasboundingbox (-1,-0.75) rectangle (3.75,4.75);
	\coordinate (u3) at (0,0);
	\coordinate (v3) at (2,0);
	\coordinate (u2) at (0,1.5);
	\coordinate (v2) at (2,1.5);
	\coordinate (u1) at (0,3);
	\coordinate (v1) at (2,3);
	\coordinate (w) at (3.5,4.5);
	\begin{pgfonlayer}{main}
		\draw (u1) node[vtx,label={[label distance=0pt]above:{\large $u_1$}}]{};
		\draw (u2) node[vtx,label={[label distance=-2pt]left:{\large $u_2$}}]{};
		\draw (u3) node[vtx,label={[label distance=0pt]below:{\large $u_3$}}]{};
		\draw (v1) node[vtx,label={[label distance=0pt]above:{\large $v_1$}}]{};
		\draw (v2) node[vtx,label={[label distance=-2pt]right:{\large $v_2$}}]{};
		\draw (v3) node[vtx,label={[label distance=0pt]below:{\large $v_3$}}]{};
		\draw (w) node[vtx,label={[label distance=0pt]left:{\large $w$}}]{};
	\end{pgfonlayer}
	\begin{pgfonlayer}{back}
		\foreach\i in {1,2,3}{
			\draw[edg,red] (u\i) -- (v\i);}
		\draw[edg,red] (v1) -- (w);
		\draw[edg] (u1) -- (u2) -- (u3);
		\draw[edg] (v1) -- (v2) -- (v3);
		\draw[edg] (u1) to[out=210, in=-210, distance=36pt] (u3);
		\draw[edg] (v1) to[out=-30, in=30, distance=36pt] (v3);
		\draw[edg] (v1) -- (v3);
	\end{pgfonlayer}
\end{tikzpicture}
\caption{$G'=C_3\square K_2$ in Case 3, $H$-edges are indicated with red.}
\end{figure}

If all $H$-edges are inside $V(G')$, then we indeed have $|E(H)|\leq 12 \leq 2\Delta$ (a planar graph on 6 vertices has at most 12 edges).
Otherwise, assume that there is an edge $e\in E(H)$ intersecting $V(G')$ in exactly one vertex. By symmetry,  $e\cap V(G')=\{v_1\}$ and let $w$ be the endpoint of $e$ different from $v_1$ (see Figure 2).
By Lemma \ref{planarF} (i) again, the three $H$-edges $e_2$, $e_3$ and $e$ again determine a $C_3\square K_2$ in $G$. If $\{u_2, u_3, w\}$ and $\{v_2, v_3, v_1\}$ were the two $C_3$-s,
then we would get a topological $K_{3,3}$ between vertex sets $\{u_2, u_3, v_1\}$ and $\{v_2, u_1, w\}$, a contradiction. Therefore, the two $C_3$-s must be
$\{u_2, u_3, v_1\}$ and $\{v_2, v_3, w\}$. Thus $w$ is connected to both $v_2$ and $v_3$ (and of course to $v_1$ by definition). This in particular implies that we
cannot have a $w'\ne w$ in the role of $w$ (namely connected to $v_1$, $v_2$ or $v_3$ by an $H$-edge) because then we would get a topological $K_{3,3}$ between vertex sets $\{v_1,v_2,v_3\}$ and $\{w,w',u_1\}$, a contradiction.

Similarly, there could be at most one $w'\in V(G)\setminus V(G')$, such that $w'$ is connected to $u_1$, $u_2$ or $u_3$ by an $H$-edge. We will show that in fact this $w'$ cannot exist assuming that $w$ exists.
If $w'$ is connected to $u_1$, then as above, by applying Lemma \ref{planarF} (i) again, $\{v_2, v_3, u_1\}$ and $\{u_2, u_3, w'\}$ are $C_3$-s. But then, $\{u_1, u_2, v_1, v_2\}$ form a $K_4$ and adding the vertex $u_3$ we get
a topological $K_5$, a contradiction.
Assume that $w'$ is connected to $u_2$ (similar for $u_3$). Applying Lemma \ref{planarF} (i) three times, we get that
$\{u_1, w', u_3\}$, $\{u_1, u_2, u_3\}$ and $\{u_2, u_3, v_1\}$ are all $C_3$-s. However, then $\{u_1, u_2, u_3, v_1, w'\}$ form a $K_5$, a contradiction. Thus $w'$ cannot exist.

Hence we can have at most 13 edges in $H$, namely the 12 potential $H$-edges induced by $V(G')$ and potentially $e=(v_1,w)$.
However, if the $H$-edge $(v_1,w)$ is present, then the edge $(u_1, u_2)$ (or the edges $(u_2, u_3)$, $(u_1, u_3)$) cannot be in $E(H)$ simultaneously.
Indeed, otherwise applying Lemma \ref{planarF} (i) again to the three $H$-edges $e_3$, $(u_1,u_2)$ and $e$, we get a topological $K_5$ again.
Therefore, $\omega(L^+(G))=|E(H)|\leq 12 \leq 2\Delta$, if $\Delta \geq 6$, as desired. \qed

\subsection{ Proof of Theorem \ref{subcubic}}

We are going to prove the theorem by induction on the number of vertices $n$.
To start the induction, note that for $n\leq 4$ we have at most 6 edges and
they could all be of different color.

Assume that the statement is true for all $n'<n$. Consider a planar
graph $G$ on $n$ vertices with maximum degree at most three and we have to find a B-coloring
of $G$ with 6 colors. We will use the colors $a,b,c,d,e,f$ and denote by $\varphi(u,v)$ the color of the edge $(u,v)$.

If $G$ does not contain a 4-cycle
then
$$q_B(G)=q(G)\leq \Delta + 1=4.$$
Thus we may assume that there is a 4-cycle $C=(v_1,v_2,v_3,v_4)$ in $G$.
Denote by $H$ the subgraph of $G$ that we get after removing $V(C)$. By induction
there is a B-coloring $\varphi$ of $E(H)$ with 6 colors. We have to extend $\varphi$ to the
rest of the edges of $G$. We distinguish the following cases. \\

{\bf Case 1:} There are 6 edges in $G|_{\{v_1,v_2,v_3,v_4\}}$.\\
In this case $\{v_1,v_2,v_3,v_4\}$ induces a $K_4$ in $G$. This is a component of $G$
so we can color this with 6 colors independently from $H$. \\

{\bf Case 2:} There are 5 edges in $G|_{\{v_1,v_2,v_3,v_4\}}$.\\
Assume wlog that the edge $(v_2,v_4)$ is missing. Then there could be an edge $e_1$ from $v_2$ to $H$
and an edge $e_2$ from $v_4$ to $H$. Let us assume first that $e_1$ and $e_2$ exist and go to different vertices $x$ and $y$
in $H$. Then we can define $\varphi(e_1)=\varphi(e_2)=a$  that is different from the at most
4 colors used on the edges incident to $x$ and $y$ in $H$. The edges $e_1$ and $e_2$ are not on a 4-cycle,
so this is allowed. Then we color the 5 edges in $G|_{\{v_1,v_2,v_3,v_4\}}$ with the remaining 5 colors.

If  $e_1$ and $e_2$ go to the same vertex $x$ in $H$, then set $\varphi(e_1)=b, \varphi(e_2)=a$
which are different from the color used on the one edge incident to $x$ in $H$ (if there is one).
Then set $\varphi(v_1,v_3)=a$  (again this is allowed, since there is no 4-cycle on the
two edges colored with $a$) and color the 4 edges of the 4-cycle $C=(v_1,v_2,v_3,v_4)$ with the remaining 4 colors
(see Figure 3). Of course, if either $e_1$ or $e_2$ (or both) is missing, the proof is even simpler. \\


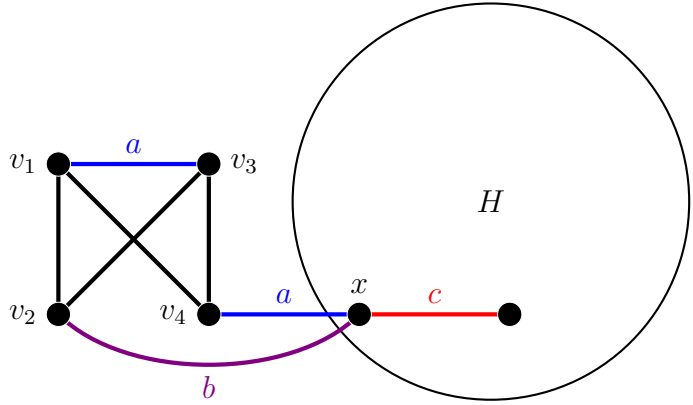
\begin{figure}[ht]
\centering
\begin{tikzpicture}[scale=1]
	\useasboundingbox (-4.75,-1.5) rectangle (4.5,4.25);
	\coordinate (v2) at (-4,0);
	\coordinate (v4) at (-2,0);
	\coordinate (x) at (0,0);
	\coordinate (y) at (2,0);
	\coordinate (h) at (1.75,1.5);
	\coordinate (v1) at (-4,2);
	\coordinate (v3) at (-2,2);
	\begin{pgfonlayer}{main}
		\draw (v1) node[vtx,label={[label distance=0pt]left:{\large $v_1$}}]{};
		\draw (v2) node[vtx,label={[label distance=0pt]left:{\large $v_2$}}]{};
		\draw (v3) node[vtx,label={[label distance=0pt]right:{\large $v_3$}}]{};
		\draw (v4) node[vtx,label={[label distance=0pt]left:{\large $v_4$}}]{};
		\draw (x) node[vtx,label={[label distance=0pt]above:{\large $x$}}]{};
		\draw (y) node[vtx] {};
		\draw (h) node {\large $H$};
	\end{pgfonlayer}
	\begin{pgfonlayer}{back}
		\draw[black,thick] (h) circle[radius=75pt] {};
		\draw[edg] (v1) -- (v2);
		\draw[edg] (v2) -- (v3);
		\draw[edg] (v1) -- (v4);
		\draw[edg] (v3) -- (v4);
		\draw[edg,red] (x) -- (y) node[midway,above] {\large $c$};
		\draw[edg,blue] (v1) -- (v3) node[midway,above] {\large $a$};
		\draw[edg,blue] (v4) -- (x) node[midway,above] {\large $a$};
		\draw[edg,violet] (v2) to[out=-45, in=-135, distance=36pt] node[below] {\large $b$} (x);
	\end{pgfonlayer}	
\end{tikzpicture}
\caption{Coloring in Case 2.}
\end{figure}

{\bf Case 3:} There are 4 edges in $G|_{\{v_1,v_2,v_3,v_4\}}$.\\
So we have the cycle $C=(v_1,v_2,v_3,v_4)$ and there are no more edges within these 4 vertices.
Then we could have one edge $e_i$ going out from each $v_i, 1\leq i \leq 4$ and not belonging to the cycle $C$. Assume wlog that
these all exist. Depending on where these 4 edges go we distinguish the following subcases. \\

{\bf Subcase 3.a:} The 4 edges $e_i$ go to 4 different vertices $x_i$.
Then we can define  $\varphi(e_1)=\varphi(e_3)=a$ that is different from the at most
4 colors used on the edges incident to $x_1$ and $x_3$ in $H$. Furthermore, set $\varphi(e_2)=\varphi(e_4)=b$ where $b$ is different from $a$ and the at most
4 colors used on the edges incident to $x_2$ and $x_4$ in $H$. We will color the remaining 4 edges
of $C$ with the remaining 4 colors $\{c,d,e,f\}$, but we have (at most) one
forbidden color for each edge:
$$\varphi(v_i,v_{i+1}) \not= \varphi(x_i,x_{i+1}), 1\leq i \leq 4$$
(counting mod 4) assuming the edges exist on the right hand side, since these would create a forbidden
4-cycle. Clearly, this can be done. For example, if
$$\varphi(x_1,x_2)=c, \varphi(x_2,x_3)=d, \varphi(x_3,x_4)=e, \varphi(x_4,x_1)=f,$$
then we can shift the coloring by one to get
$$\varphi(v_1,v_2)=d, \varphi(v_2,v_3)=e, \varphi(v_3,v_4)=f, \varphi(v_4,v_1)=c.$$
This results in a B-coloring. \\

{\bf Subcase 3.b:} $x_1=x_3$ and $x_2=x_4$.
Note that there cannot be an edge between $x_1$ and $x_2$ in $H$, since that would create a $K_{3,3}$
between the two sets $\{v_2,v_4,x_1\}$ and $\{v_1,v_3,x_2\}$.
Set $\varphi(e_1)=a, \varphi(e_3)=b$ where $a$ and $b$ are different from the 2 colors used on the one edge incident to $x_1$ in $H$ and the one edge incident to $x_2$ in $H$  (if these exist).
Furthermore, set $\varphi(e_2)=a, \varphi(e_4)=b$.
Then we can color the remaining 4 edges
of $C$ with the remaining 4 colors $\{c,d,e,f\}$, because we have no restrictions in this subcase. \\

{\bf Subcase 3.c:} $x_1=x_2$ and $x_3=x_4$.
If the edge $(x_1,x_3)$ exists, then the 4-cycle together with these two points form a 3-regular component
which clearly can be colored with 6 colors. Otherwise,
we color $e_1$ and $e_3$ with a color, say $a$, that is different
from the 2 colors used on the one edge incident to $x_1$ in $H$ and the one edge incident to $x_3$ in $H$.
We color $e_2$ and $e_4$ with a color that is different from $a$ and the at most
the two colors used on the one edge incident to $x_1$ in $H$ and the one edge incident to $x_3$ in $H$. We can color the remaining 4 edges
of $C$ with the remaining 4 colors $\{c,d,e,f\}$, because we have no restrictions again. \\

{\bf Subcase 3.d:} $x_1=x_3$, but $x_2\not=x_4$  ($x_2=x_4$, but $x_1\not=x_3$ is symmetric).
Note that one of the edges $(x_1,x_2)$ and $(x_1,x_4)$ must be missing, otherwise
the degree of $x_1$ would be at least 4. Wlog assume that $(x_1,x_4)$ is missing. We may also assume that
$(x_1,x_2)$ is an edge otherwise the proof is even easier.
Set $\varphi(e_1)=a, \varphi(e_3)=b$ where $a$ and $b$ are different from the 4 colors used on the edge $(x_1,x_2)$, on the one possible additional edge incident to $x_2$ in $H$ and the two edges incident to $x_4$ in $H$  (if these exist).
Set $\varphi(e_4)=a,\varphi(e_2)=c$ where $c$ is different from $a,b,\varphi(x_1,x_2)$ and from the color of the possible further edge of $H$ incident to $x_2$.
Set $\varphi(v_1,v_4)=c$ and color the remaining 3 edges of $C$ with the remaining 3 colors $\{d,e,f\}$; with the restriction
$$\varphi(v_1,v_2),\varphi(v_2,v_3) \not= \varphi(x_1,x_2).$$ \\

{\bf Subcase 3.e:} $x_3=x_4$, but $x_1\not=x_2$  ($x_1=x_2$, but $x_3\not=x_4$ is symmetric).
Note that out of the edges $(x_1,x_3)$ and $(x_2,x_3)$ one of them must be missing, otherwise
the degree of $x_3$ would be at least 4. Wlog assume that $(x_1,x_3)$ is missing. We may also assume that
$(x_2,x_3)$ is an edge otherwise the proof is even easier.
We select two colors, say $a$ and $b$,
which are different from the 4 colors used on the edge $(x_2,x_3)$, on the one possible additional edge incident to $x_2$ in $H$ and the two edges incident to $x_1$ in $H$  (if these exist).
Set  $\varphi(e_1)=\varphi(e_3)=a$  and $\varphi(e_2)=\varphi(_4)=b$.
Color the remaining 4 edges
of $C$ with the remaining 4 colors $\{c,d,e,f\}$; we have two restrictions
$$\varphi(v_1,v_2) \not= \varphi(x_1,x_2), \varphi(v_2,v_3)\not=\varphi(x_2,x_3).$$

{\bf Subcase 3.f:} Three of the $x_i$-s are the same and the 4-th is different, say $x_1=x_3=x_4$, but $x_2\not=x_1$ (otherwise symmetric).
Note that the edge $(x_1,x_2)$ must be missing, otherwise
the degree of $x_1$ would be at least 4.
We select a color, say $a$,
which is different from the 2 colors used on the two edges incident to $x_2$ in $H$  (if these exist).
We color $e_1$ and $e_2$ with color $a$, $e_3$ with color $b$ and $e_4$ with color $c$.
We color the remaining 4 edges
of the 4-cycle with the 4 colors $\{c,d,e,f\}$; we have two restrictions
$$\varphi(v_1,v_4) \not= c, \varphi(v_3,v_4)\not=c$$
(otherwise the coloring would not be proper). Again this clearly can be done. There are no more cases
and this finishes the proof of the theorem. \qed

\subsection{Proof of Theorem \ref{outerplanar}}
The two exceptions in (i) are obvious. We are going to prove the non-exceptional part of (i) by induction on $n=|V(G)|$. Assume that $G$ is a non-exceptional (i.e not $K_4-e$ or $C_4$),
$2$-connected, outerplanar graph on $n$ vertices with maximum degree at most $\Delta$. Let $G$ consist of the outer
cycle $C = (v_1, v_2, \ldots , v_n)$ and some chords of $C$. We may assume that $\Delta(G)\ge 3$, i.e. there are some chords in $C$ otherwise $G=C$ and apart from the excluded $C_4$, we have $q_B(C_i)\le \Delta(C_i)+1$, finishing the proof.

We can start the induction with $n=3$. Then $G=C_3$ and (i) clearly holds.  Assume that (i) holds for all $3<n'<n,n'\ne 4$.  We call a face a {\em boundary face} if it consists of one chord
and a subpath of the outer cycle $C$ (i.e. there are no additional chords
within the face). Consider a boundary face with chord $(v_i,v_j)$ (where $v_i,v_j$ are not consecutive on $C$).
Denote by $H$ the subgraph we get by removing from $G$ the inner vertices of the path $P$ of $C$ connecting $v_i$ and $v_j$ on the boundary face.
Clearly, $H$ is an outerplanar graph.

If $|V(H)|=4$, i.e. $H$ is a $C_4$ or a $K_4-e$, then their B-colorings (with 4 or 5 colors, respectively) can easily be extended to a B-coloring of $G$ (with 4 or 5 colors, respectively), proving (i).
Otherwise, by the inductive hypothesis, $H$ has a B-coloring $\varphi$  with $\Delta+1$ colors. The strategy of the proof is to extend this coloring $\varphi$ to $P$ without using
any additional colors.

Since  $d_H(v_i),d_H(v_j)\le \Delta-1$, we have at most $\Delta-1$ forbidden colors to color edges of $P$, at both $v_i$ and $v_j$.
In addition, if this boundary face
is a 3-face (so $j=i+2$) and there is an {\em inner}
triangle (facing the other direction away from the boundary face) sitting on the $(v_i,v_j)$ edge,
say $(v_i, v_j, v_l)$, then this gives a forbidden color. Indeed, we have
the restrictions (see Figure 4):
$$\varphi(v_i,v_{i+1}) \not= \varphi(v_{i+2},v_l), \varphi(v_{i+1},v_{i+2})\not=\varphi(v_i,v_l),$$
since otherwise we get a 4-cycle $(v_i, v_{i+1}, v_{i+2}, v_l)$ with a repeated color.


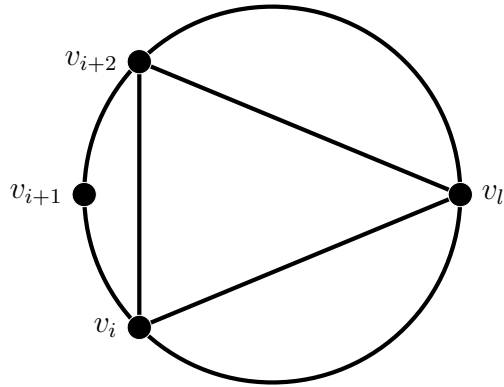
\begin{figure}[ht]
\centering
\begin{tikzpicture}[scale=1]
	\useasboundingbox (-3.5,-2.5) rectangle (3.25,2.5);
	\coordinate (vi2) at (135:2.5);
	\coordinate (vi1) at (180:2.5);
	\coordinate (vi0) at (-135:2.5);
	\coordinate (vl) at (0:2.5);
	\begin{pgfonlayer}{main}
		\draw (vi2) node[vtx,label={[label distance=0pt]left:{\large $v_{i+2}$}}]{};
		\draw (vi1) node[vtx,label={[label distance=0pt]left:{\large $v_{i+1}$}}]{};
		\draw (vi0) node[vtx,label={[label distance=0pt]left:{\large $v_{i}$}}]{};
		\draw (vl) node[vtx,label={[label distance=0pt]right:{\large $v_{l}$}}]{};
	\end{pgfonlayer}
	\begin{pgfonlayer}{back}
		\draw[edg] (0,0) circle[radius=2.5] {};
		\draw[edg] (vl) -- (vi0) -- (vi2) -- (vl);
	\end{pgfonlayer}	
\end{tikzpicture}
\caption{The extra restriction from the inner triangle.}
\end{figure}

Furthermore, if the boundary face is a 4-face (so $j=i+3$) we have one more restriction
$$\varphi(v_{i+1},v_{i+2}) \not= \varphi(v_i,v_j).$$

If the boundary face is a cycle of length at least four, we can extend $\varphi$ to edges of $P$ as follows.
Let us select two colors $a$ and $b$, such that $a\ne b$ and $a,b$ are not used in the coloring of $H$ at $v_i,v_j$, respectively.
The edges of $P$ can be colored with $a,b$ and an arbitrary third color different from $\varphi(v_i,v_j)$ (needed only for length four), maintaining a B-coloring of $G$.  Thus we may assume that every boundary face is a 3-face. This argument works also on a boundary 3-face if either $d_G(v_i)$ or $d_G(v_{i+2})$ is smaller than $\Delta$. Indeed, in this case we color first the edge at the endpoint of $P$ where there is only one free color, say $a$ (because its degree is $\Delta$ and the inner triangle gives another restriction). Then for the edge of $P$ at the other endpoint (with degree smaller than $\Delta$) there is still an available color different from $a$.

Thus we may assume that the endpoints of the chords at {\em all} boundary 3-faces are of degree $\Delta$. We may also assume that all boundary 3-faces have an inner triangle at the two endpoints $v_i$ and $v_{i+2}$, otherwise we can extend the coloring with the above technique.  Note that an inner triangle is really an inner triangle; we cannot have
for example $(v_{i-1},v_i,v_{i+2})$ as the inner triangle because then $v_i$
cannot have any additional chords since they would cross the third side of the
inner triangle. Thus $v_i$ has degree $\Delta=3$, implying that $v_{i+2}$ also has degree 3, so $G=K_4-e$, a contradiction. Also, since we have inner triangles, $\Delta \ge 4$.

We call two boundary triangles {\em neighbors}, if they share an endpoint.
We will need the following.

\begin{claim}\label{inner}
There are two neighboring boundary triangles that have the same inner triangle.
\end{claim}

\noindent {\bf Proof of Claim \ref{inner}:}
We define the {\em distance} of a boundary triangle from the third point of its
inner triangle; out of the two subpaths of $C$ from the endpoints to the third point of the
inner triangle we take the length (counting edges) of the shorter one.
Consider the boundary triangle with the minimum distance (at least 2 by the above), say $v_i$ and $v_{i+2}$
and the third vertex of the inner triangle is $v_l$. Say the distance is given by the $(v_i,v_l)$ subpath $P_d$ on $C$.
Consider the region $R$ bounded by the chord $(v_i, v_l)$ and the subpath $P_d$. This region
also must contain a boundary triangle. If this is a boundary triangle itself, then we are done.
Otherwise, $R$ properly contains a boundary triangle along with its inner triangle.
Then the distance of this boundary triangle to the third point of its inner triangle would be less, a contradiction, proving the claim. \qed

From the claim we immediately get a contradiction if $\Delta \ge 5$.  Thus we have $\Delta =4$ and  we proceed as follows.
We will color with colors $\{a,b,c,d,e\}$.
Assume that the two neighboring boundary triangles from the claim
are $(v_1, v_2, v_3)$ and $(v_3,v_4, v_5)$ along with the common inner triangle $(v_1,v_3,v_5)$.
If $V(G)=\{v_1,v_2,v_3,v_4,v_5\}$ then $G$ has a B-coloring with $5$ colors: the inner triangle is colored with 1,2,3 the edges $(v_1,v_2),(v_3,v_4)$ are colored with 4 and the edges  $(v_2,v_3),(v_4,v_5)$ are colored with 5.  Otherwise $V(G)\setminus \{v_2,v_3,v_4\}$ spans a $2$-connected outerplanar $H'$ with Hamiltonian cycle $(v_1,v_5,v_6,\dots,v_n)$.

By induction on $n$, $H'$ has a B-coloring $\varphi'$  with at most $\Delta+1= 5$ colors. This is also true if $H'$ is exceptional, i.e. $n=7$ and $(v_1,v_5,v_6,v_7)$ is a $C_4$. (The case $H'=K_4-e$ is not possible since $\Delta=4$.)   We have to extend $\varphi'$ to the rest of $G$ (6 additional edges).  Assume wlog that $\varphi'(v_1,v_5)=a,\varphi'(v_1,v_n)=b$. \\

{\bf Case 1:} $\varphi'(v_5,v_6)=c$. We can extend $\varphi'$ as follows.
$$\varphi'(v_1,v_2)=\varphi'(v_3,v_4)=c,\varphi'(v_2,v_3)=\varphi'(v_4,v_5)=b,$$
$$\varphi'(v_1,v_3)=d, \varphi'(v_3,v_5)=e.$$ \\

{\bf Case 2:} $\varphi'(v_5,v_6)=b$. Since $\varphi'$ is a B-coloring on $H'$, we have $v_6\ne v_n$ and $(v_6,v_n)\notin E(H')$. Therefore we can change the color $b$ on $(v_1,v_n)$ or on $(v_5,v_6)$ to one of the colors $\{c,d,e\}$ so that $\varphi'$ is still a B-coloring on $H'$, unless the edges of $H'$  incident to $v_n$ and to $v_6$ are both colored with colors $\{b,c,d,e\}$. However, in this case the recoloring
$$\varphi'(v_1,v_5)=c, \varphi'(v_5,v_6)=a$$
gives a B-coloring on $H'$ (using that $(v_n,v_6)\notin E(H')$). This leads back to Case 1 and finishes the proof of (i).
\medskip

To prove (ii), i.e. $\omega(L^+(G))\le \Delta$, we use Lemma \ref{outerplanarF}. Let $K$ be a complete subgraph in $L^+(G)$ and let $H$ be the subgraph of $G$ determined by the edge set corresponding to the vertices of $K$. Let $S$ be the largest complete subgraph of $K$ formed by the edges of $F(G)$. By Lemma \ref{outerplanarF} (ii), $|S|\le 2$. \\

{\bf Case 1:} $|S|=1$. Here a single edge $e=(u_1,v_1)$ corresponds to $S$ and all further edges of $H$ must meet $e$. However, since $H$ has no two independent edges, $H$ must be a star, or a triangle.  Thus $\omega(L^+(G))=|E(H)|\le \max\{ \Delta(H),3\} \le \Delta$. \\

{\bf Case 2:} $|S|=2$. We will show that in this case $|E(H)|\leq 5 \leq \Delta$.
Here edges $e_1=(u_1,v_1),e_2=(u_2,v_2)$ correspond to $S$ and $G$  has a 4-cycle $C=(u_1,u_2,v_2,v_1)$ with edges $e_1,e_2$ on it.
If all $H$-edges are inside $V(C)$, then we indeed have $|E(H)|\leq 5 \leq \Delta$. Thus we may assume that we have additional $H$-edges.
Clearly no edge $e \in E(H)$  exists with $|e\cap C|=0$ otherwise $|S|>2$, contradicting Lemma \ref{outerplanarF}. Since outerplanar graphs have no $K_4$ subgraphs, we may assume that the edge $(u_1,v_2)$ is not in $G$. Moreover, no edge $f\in E(H)$ can intersect $C$ in exactly the vertex $u_1$ otherwise the edge $f$ and $(u_2,v_2)$ cannot be in a 4-cycle without forming a $K_{2,3}$ subgraph which is also impossible in an outerplanar graph. The same argument works for $v_2$, thus all edges of $H$ meet $C$ in either $v_1$ on in $u_2$. However we claim that there are no two edges $f,g\in E(H)$ meeting $C$ in $v_1,u_2$, respectively. Indeed, $f,g$ must be independent, otherwise we have a $K_{2,3}$. But $f,g$ can be in a 4-cycle of $G$ only if $(u_2,v_1)\in E(G)$ and the other endpoints of $f,g$ are adjacent in $G$. But then this edge with $f,g$ give a path between $u_2,v_1$ extending the paths $(u_2,u_1,v_1)$ and $(u_2,v_2,v_1)$  to a topological $K_{2,3}$.
Thus, apart from at most two edges all edges of $H$ are incident $v_1$ or to $u_2$, assume wlog that to $v_1$.

Let $f=(v_1,z)$ be an edge of $H$ intersecting $C$ in only $v_1$. The edges  $f,e_2\in E(H)$ are in a 4-cycle of $G$. The edge $(u_2,z)\notin E(G)$ since there is no $K_{2,3}$ in $G$. Thus, for all such $z$ the 4-cycle must use the edges $e_2,(v_2,z),f,(v_1,u_2)$. But then there could be only one such $z$ otherwise we have a $K_{2,3}$ in $G$, giving us at most 6 $H$-edges.
Furthermore, the edge $(u_1,u_2)$ cannot belong to $E(H)$ assuming that this $z$ exists, since otherwise the edges $(u_1,u_2)$ and $f$ are in a 4-cycle and then we either get a $K_4$ or a $K_{2,3}$.
Thus indeed $\omega(L^+(G))=|E(H)|\le 5\le \Delta$ if $\Delta \ge 5$. This concludes the proof of Theorem \ref{outerplanar}. \qed

\smallskip

\noindent{\bf Acknowledgement. } The authors appreciate the remarks of the referees that improved the presentation.

\end{document}